\documentclass[11pt]{amsart}
\usepackage{calc,amssymb,amsthm,amsmath,fullpage}
\RequirePackage[dvipsnames,usenames, table]{xcolor}
\usepackage{hyperref}
\usepackage[capitalize]{cleveref}
\hypersetup{
  bookmarks,
  bookmarksdepth=3,
  bookmarksopen,
  bookmarksnumbered,
  pdfstartview=FitH,
  colorlinks,backref,hyperindex,
  linkcolor=Sepia,
  anchorcolor=BurntOrange,
  citecolor=MidnightBlue,
  citecolor=OliveGreen,
  filecolor=BlueViolet,
  menucolor=Yellow,
  urlcolor=OliveGreen
}
\usepackage{xypic}
\usepackage{alltt}
\usepackage{multicol}  
\usepackage{xspace}
\usepackage{rotating}
\theoremstyle{definition}
\newtheorem{example}{Example}[section]
\newtheorem{remark}{Remark}[subsection]

\newcommand{\bA}{\mathbb{A}}
\newcommand{\bF}{\mathbb{F}}

\usepackage{verbatim}
\usepackage{enumerate}
 
\begin{document}
\title{{RandomPoints} package for \emph{Macaulay2}}
\author{Sankhaneel Bisui}
\email{sbisui@tulane.edu}
\address{Department of Mathematics, Tulane University, New Orleans, LA 70118}
\author{Zhan Jiang}
\email{zoeng@umich.edu}
\address{Department of Mathematics, University of Michigan, Ann Arbor, MI 48109}
\author{Sarasij Maitra}
\email{sm3vg@virginia.edu}
\address{Department of Mathematics, University of Virginia, Charlottesville, VA 22904}
\author{Th\'ai Th\`anh Nguy$\tilde{\text{\^e}}$n}
\email{tnguyen11@tulane.edu}
\address{Department of Mathematics, Tulane University, New Orleans, LA 70118}
\author{Karl Schwede}
\email{schwede@math.utah.edu}
\address{Department of Mathematics, University of Utah, 155 S 1400 E Room 233, Salt Lake City, UT, 84112}

\thanks{Schwede was supported by NSF Grants \#1801849, FRG \#1952522 and a Fellowship from the Simons Foundation.}

\begin{abstract}
  We present {\tt RandomPoints}, a package in \emph{Macaulay2} designed mainly to identify rational and geometric points in a variety over a finite field.  We provide tools to estimate the dimension of a variety.  We also present methods to obtain non-vanishing minors of a given size in a given matrix, by evaluating the matrix at a point.  
\end{abstract}


\maketitle

\section{Introduction}
    Let $I$ be an ideal in a polynomial ring $k[x_1,\dots, x_n]$ over a finite field $k$. Let $X:=V(I)$ denote the corresponding set of rational points in affine $n$-space. Finding one such rational point or geometric point (geometric meaning a point over some finite field extension), in an algorithmically efficient manner is our primary motivation for this package. The authors of the package are Sankhaneel Bisui, Zhan Jiang, Sarasij Maitra, Th\'ai Th\`anh Nguy$\tilde{\text{\^e}}$n, Frank-Olaf Schreyer, and Karl Schwede.

    There is an existing package called {\tt RationalPoints} \cite{RationalPointsSource}, which we took inspiration from, which aims to find \emph{all} the rational points of a variety; our aim here is to find one or more points quickly, even if it is not rational.  We note also that the package {\tt Cremona} \cite{CremonaSource} can find rational points on projective varieties, as can built-in the function {\tt randomKRationalPoint} \cite{M2}, our methods frequently appear to be faster and apply in the affine setting as well.

    We develop functions that apply different strategies to generate random rational and geometric points on the given variety. We also provide functions that will expedite the process of determining properties of the singular locus of $X$.

    We provide the following core functions:
    \begin{itemize}    
    \item {\tt randomPoints}:  This tries to find a point in the vanishing set of an ideal. (\Cref{randomPoints})
    \item {\tt dimViaBezout}:  This tries to compute the dimension of an algebraic set by intersecting with hyperplanes.  (\Cref{subsec.DimViaBezout})
    \item {\tt projectionToHypersurface} and {\tt genericProjection}: These functions provide customizable projection. (\Cref{projectionfunctions}) 	
    \item {\tt findANonZeroMinor} and {\tt extendIdealByNonZeroMinor}:  The first function finds a submatrix of a given matrix that is nonsingular at a point of a given ideal.  The second, adds said submatrix to an ideal, which is useful for computing partial Jacobian ideals.  (\Cref{findANonZeroMinor})          
    \end{itemize}

    All polynomial rings considered here will be over finite fields. In the subsequent sections, we explain  the core and helper functions and describe the strategies that we have implemented.

   \vspace{1em}
    \noindent \textbf{Acknowledgements:} The authors would like to thank David Eisenbud and Mike Stillman for useful conversations and comments on the development of this package.  The authors began work on this package at the virtual Cleveland 2020 Macaulay2 workshop. The authors are also grateful to the reviewers for suggesting and providing preliminary codes to speed up computations, thereby improving the efficacy of the package substantially.  Frank-Olaf Schreyer is also an author on this package, as he provided code related to the {\tt MultiplicationTable} decomposition strategy and suggested using a probabilistic approach to compute dimension.

\section{Our primary purpose: {\tt randomPoints}}\label{randomPoints}
    We start with the core function in this package:  {\tt randomPoints} is a function to find rational or geometric points in a variety. The typical usages are as follows: 

    \vspace{0.5em}
    -- {\tt randomPoints(I)}, 

    -- {\tt randomPoints(n, I)} 

    \vspace{0.5em}
    \noindent where $n$ is a positive integer denoting the number of desired points, and 
    $I$ is an ideal inside a polynomial ring.  If {$ n$} is omitted, then it is assumed to be 1.

    \subsection{Options}\label{strategydetails}

    The user may also choose to provide some additional information depending on the context which may help in faster computations, or whether a point is found at all.

    \begin{description}
        \setlength{\itemsep}{5pt}
    \item[\tt Strategy => $\bullet$]   Here the $\bullet$ can be {\tt Default, BruteForce} or {\tt LinearIntersection} \\ .

    \begin{itemize}
    \item {\tt Default} performs a sequence of the different strategies below, aimed at finding a point quickly.  It begins begins with brute force and moves to linear intersections with particularly simple linear forms, and gradually ramps up the randomness.  
    \item {\tt BruteForce} simply tries random points and sees if they are on the variety.
           
    \item {\tt LinearIntersection} intersects with a random linear space.  
    \end{itemize}

    Notice that the speed, or success, varies depending on the strategy (see also \Cref{helper}).

    \begin{example}\label{BruteForce}
        Consider the following example.
        ~~        
        {\small\color{blue}
    \begin{quote}
\begin{verbatim}
i2 : R = ZZ/101[x, y, z];

i3 : J = ideal(x^3 + y^2 + 1, z^3 - x^2 - y^2 + 2);

o3 : Ideal of R

i4 : time randomPoints(J,Strategy=>BruteForce, PointCheckAttempts=>10)
        -- used 0.00186098 seconds        

o4 = {}

o4 : List

i5 : time randomPoints(J)
        -- used 0.0205099 seconds

o5 = {{-1, 0, -1}}

o5 : List

i6 : time randomPoints(J, Strategy=>LinearIntersection)
        -- used 0.0334881 seconds                

o6 = {{0, 10, 48}}
\end{verbatim}
    \end{quote}
        }
    \end{example}%
    \vspace{-1em}

    \item[\tt ExtendField => Boolean] 

    Intersection with a general linear space will naturally find scheme theoretic points that are not rational over the base field.  
    Setting {\tt ExtendField => true} will tell the function that such points are valid.  Setting {\tt ExtendField => false} will tell the function ignore such points.  
    In fact, {\tt ExtendField => true} will also tell Macaulay2 to use linear spaces defined over a field extension, which can improve randomness properties.
    This sometimes can slow computation, and other times can substantially speed it up when the variety has few rational points.  For some applications, points over extended fields may also have better randomness properties.

    \item[\tt DecompositionStrategy]  Within the {\tt LinearIntersection} strategy, one can also specify the option {\tt DecompositionStrategy}.  Valid values are {\tt Decompose} and {\tt MultiplicationTable}, the latter of which is currently only implemented for homogeneous ideals.  The point is, after intersecting the linear space and obtaining an ideal defining a set of (possibly thickened) points, we need to find the minimal associated primes.  By default we use Macaulay2's built-in {\tt decompose} command.  We also have implemented a {\tt MultiplicationTable} algorithm, as provided by Frank-Olaf Schreyer, which utilizes the action of a variable on the residue fields of these points computed in more than one way.  This method is frequently faster for rings with smaller numbers of variables.
    
    The {\tt Default} strategy switches back and forth between {\tt Decompose} and {\tt MultiplicationTable} for homogeneous ideals (starting with one the function thinks will be fastest).  Setting this to {\tt Decompose} in the default strategy will force only {\tt Decompose} to be used; setting it to {\tt MultiplicationTable} will force only {\tt MultiplicationTable} to be used (if the ideal is homogeneous).      

    \item[\tt Homogeneous]  Setting this true specifies that the origin (corresponding to the irrelevant ideal) is not a valid point.

    \item[\tt Replacement => Monomial, Binomial, Trinomial, Full]
        When intersecting with a random linear space, it is frequently much faster to use a linear space defined by relatively sparse equations (ie, equations that do not involve all variables).   Specifying {\tt Monomial} will mean linear forms such as $ax + b$ are used (for constants $a$ and $b$), involving only one variable.  Binomial means forms like $ax + by + c$, using two variables.  Trinomial means forms like $ax + by + cz + d$.  {\tt Full} means all variables will have coefficients.  

    \item[\tt DimensionFunction => Function]
        Our current implementation does not need to know the dimension of $V(I)$.  However, there are places where we try to verify the dimension of an ideal before we decompose the ideal.  You can pass this function {\tt dim} (the default) or our probabilistic {\tt dimViaBezout} or any other dimension function you might prefer.

    \item[\tt PointCheckAttempts => ZZ]

    When calling {\tt randomPoints} with a {\tt BruteForce} strategy, this denotes the number of trials for brute force point checking.  It also controls how many linear spaces to simultaneously study in the {\tt LinearIntersection} strategy.
    
    \begin{example}
        We re-compute \cref{BruteForce} this time specifying more attempts.
    {{\small\color{blue}
    ~~
    \begin{quote}
\begin{verbatim}    
i7 : time randomPoints(J, Strategy => BruteForce, PointCheckAttempts => 10000)
-- used 1.16294 seconds

o7 = {{-43, 25, 29}}
    \end{verbatim}%
\end{quote}\vspace{-1em}%
}}%
    \end{example}%
    
    \item[\tt NumThreadsToUse => ZZ]

    When calling {\tt randomPoints} with a {\tt BruteForce} strategy, this denotes the number of threads to use in brute force point checking.
    \end{description}

    \subsection{Comments on performance and implementation}
    \label{subsec.CommentOnPerformanceAndImplementation}

    When working over very small fields, especially with hypersurfaces, frequently {\tt BruteForce} is most efficient.  This is not surprising as there may not be many points to check.  However, if the field size is larger, {\tt BruteForce} will perform poorly.  Even for some simple examples, it could not provide any rational points if the number of trials is not large enough.  Other strategies work differently on different examples, and the same strategy can sometimes work very quickly even if it typically works very slowly.  

    The current version of the {\tt LinearIntersection} strategy no longer computes the dimension of the algebraic set.  Instead, it first finds a point defined by linear equations.  If the point is on the algebraic set, we are done.  If not, we throw away one of the forms and so now have a line and we see if this line intersects our algebraic set.  We continue in this way until we find a point.  This appears to avoid a number of bottlenecks in our previous implementation since Macaulay2 is relatively fast at identifying when a linear space and a variety do \emph{not} intersect.  
    

    \begin{example}
        We begin with an example over a small field.
    {{\small\color{blue}
    ~~
    \begin{quote}
\begin{verbatim}
i2 : R = ZZ/7[x_1..x_10];

i3 : I = ideal(random(2, R), random(3, R));

o3 : Ideal of R

i4 : time randomPoints(I, Strategy => BruteForce, PointCheckAttempts => 20000)
        -- used 0.00311884 seconds        

o4 = {{-1, -1, 0, 2, 2, -2, -2, -3, -3, -3}}

o4 : List

i5 : time randomPoints(I, Strategy => Default)
        -- used 0.081349 seconds          

o5 = {{3, 0, 3, 3, 2, -2, 1, -1, 3, 1}}
\end{verbatim}
\end{quote}\vspace{-1em}
    }}
    \end{example}    

\begin{example}
    Now we work over a larger field.
    {{\small\color{blue}
    \begin{quote}
\begin{verbatim}
i6 : S = ZZ/211[x_1..x_10];

i7 : J = ideal(random(2, S), random(3, S));

o7 : Ideal of S

i8 : time randomPoints(J, Strategy => BruteForce, PointCheckAttempts => 2000000)
        -- used 17.7988 seconds
                   
o8 = {{15, 67, -27, -103, 56, 66, -23, 28, -50, 13}}

o8 : List

i9 : time randomPoints(J, Strategy => Default)
        -- used 0.0864013 seconds         

o9 = {{0, 0, 0, 0, 34, 76, 51, 0, 1, 0}}
    \end{verbatim}
\end{quote}\vspace{-2em}
        }}
    \end{example}

\begin{example}
        Finally, we can allow our functions to extend our field.
        {{\small\color{blue}
        \begin{quote}
            {\small
\begin{verbatim}
i11 : time randomPoints(J, Strategy => Default, ExtendField => true)
            -- used 0.144332 seconds              
       
              3      2                      3      2                        
o11 = {{0, - a  + 62a  - 47a - 76, 0, 0, 13a  - 18a  + 63a - 31, 0, 0, 
             3      2                 3      2
        - 20a  - 82a  + 35a - 19,  55a  - 64a  - 8a - 50, 1}}

i12 : coefficientRing ring first first o11

o12 = GF 1982119441       

i13 : log_211 1982119441

o13 = 4
\end{verbatim}
            }
    \end{quote}
        }}
In this case, we found a degree 4 point.
    \end{example}

\begin{remark}[Comments on the probability of finding a point]
	In the case of an absolutely irreducible hypersurface in $\bA^n_{\bF_q}$ (defined by $f$ say), there is significant discussion in the literature estimating lower bounds of number of rational points (see for instance, \cite{RationalKPointsSource,lang1954number,ghorpade2002number,cafure2006improved}) all of which point to the fact that there is ``good probability" of finding a rational point in this case when we intersect with a random line. Heuristically, we can make the following rough estimation.  We expect that each equation $f = \lambda$ for $\lambda \in \bF_q$ has approximately the same number of solutions.  Since each point on $\bF_q^n$ solves exactly one of these equations, we expect that $f = 0$ has approximately $q^{n-1}$ solutions, or in other words, our hypersurface has $q^{n-1}$-points.  
       Now, a random line $L$ has $q$ points. We want to find the probability that one of these points is  rational for $V(f)$.  We would expect that if these points are randomly distributed, then the probability that our line contains one of those points $1-(1-{1 \over q})^q$ which tends to $1-e^{-1} \approx 0.63$ for $q$ large.   Alternately, one can use the proof of \cite[Proposition 2.12]{bothmer2005quick} for a more precise statement.  For each point of $L$, we see that the probability that  the chosen point does not lie in the intersection, $L\cap V(f)$, is $1-\frac{1}{q}$. We then exhaust this search over all the points on $L$ to get the probability that there is indeed a successful intersection is $1-(1-\frac{1}{q})^{q}$. As $q$ gets larger, this value tends to $1-e^{-1}\approx 0.63$.

       Of course, there are schemes over $\bF_q$ with no rational points at all, even for plane curves.
\end{remark}

\begin{remark}[Projecting to a hypersurface first]
    \label{rem.ProjectingToAHypersurfaceFirst}
    Suppose $X \subseteq \bA^n$ is an algebraic set.  In a number of existing algorithms, one first does a generic (or even not very generic) projection $h : \bA^n \to \bA^m$ and so that $h(X)$ is a hypersurface (at least set theoretically).  Then one finds a point $x \in h(X)$ (say as above), and computes $h^{-1}(\{ x \})$ which is a linear space in $\bA^n$ that typically intersects $X$ in a rational point.   For example, this is done in {\tt randomKRationalPoint} in core Macaulay2.  Note that projecting to a hypersurface \emph{still is performing an intersection with a linear space}, since $h^{-1}(\{ x \})$ is linear, but it tries to choose the linear space intelligently.

    However, in our experience, doing this generic projection first yields slower results. First, one has to compute the dimension.  There are also numerous cases where computing this hypersurface $h(X)$ can be quite slow.  This particularly appears in cases when one is computing successive minors to identify the locus where some variety is nonsingular.  

    On the other hand, instead of using a truly random linear space to intersect with, in the default strategy, we normally choose a linear space whose defining equations are particularly non-complicated (i.e., with many linear or binomial terms).  Such simple linear spaces are the ones implicitly considered in {\tt randomKRationalPoint} for instance since that generic projection is so simple.  In practice, our approach seems to give at least as good performance compared to projecting to a hypersurface, without the chance of the code getting hung up on the generic projection or dimension computations.  We also do successive intersections in a way that avoids computing the dimension as described above in \Cref{subsec.CommentOnPerformanceAndImplementation}.
\end{remark}

\section{Useful functions {\tt dimViaBezout} and {\tt randomCoordinateChange}}{\label{helper}}

\subsection{\tt dimViaBezout:} 
\label{subsec.DimViaBezout}
We thank the Frank-Olaf Scheyer for pointing out that in most of the computations, computing the codimension of the given ideal, is a significant bottleneck. While we have avoided most dimension computations in our current implementation, we have also implemented a probibalistic dimension computation of $V(I)$.
This function takes as input an ideal $I$ in a polynomial ring over a field and intersects $V(I)$ with successively higher dimensional random linear spaces until there is an intersection.  For example, if $V(I)$ intersect a random line has a point, then we expect that $V(I)$ contains a hypersurface.  If there was no intersection, this function tries a $2$-dimensional linear space, and so on.  This can speed up a number of computations.  The function also takes in optional inputs as described below:

\begin{itemize}
	\item {\tt DimensionIntersectionAttempts}: 
        Our function actually estimates dimension several times and then averages the result (rounding down) since we tend to overestimate the dimension by the above.  By default it does this 3 times unless the {\tt Homogeneous} flag is set, in which case it is done 5 times.
	\item {\tt MinimumFieldSize}:
	If the ambient field is smaller than this integer value, it will automatically be replaced with an extension field.  For instance, there are relatively few linear spaces over a field of characteristic $2$, and this can cause incorrect results to be returned to the user. The user may set the {\tt MinimumFieldSize} to ensure that the field being worked over is big enough.  If this is not set, the program tries to choose a reasonable minimum field size based on the ambient ring.
    \item {\tt Homogeneous}: 
    If the ideal is homogeneous, we can use homogeneous linear spaces to compute dimension. Sometimes this is faster and other times slower.
\end{itemize}

\begin{example} We illustrate the speed difference in this example.
	{\color{blue}{\small
			\begin{quote}
				\begin{verbatim}
i2 : S = ZZ/101[y_0..y_9];

i3 : I=ideal random(S^1,S^{-2,-2,-2,-3})+(ideal random(2,S))^2;

o3 : Ideal of S

i4 : time dimViaBezout I
        -- used 0.837359 seconds        

o4 = 5

i5 : time dim I
        -- used 36.8496 seconds           

o5 = 5

i6 : time dimViaBezout(I, DimensionIntersectionAttempts=>1)
        -- used 0.280803 seconds  

o6 = 5
				\end{verbatim}
			\end{quote}
	}}
As you can see doing a single intersection attempt is about three times faster, and it usually gives the right answer (far more than 99\% of the time in this particular example, but in others doing the computation in triplicate avoids returning incorrect answers).
\end{example}

\subsection{\tt randomCoordinateChange:} 
\label{subSec.randomCoordinateChange}

This function takes a polynomial ring as an input and outputs a coordinate change map,
i.e. given a polynomial ring, this will produce a linear automorphism of the ring.  This function checks whether the map is an isomorphism by computing the Jacobian.

In some applications, a full random change of coordinates is not desired, as it might cause code to run very slowly.  A binomial change of coordinates might be preferred, or we might only replace some monomials by other monomials.  
This is controlled with the following options.

\begin{itemize}
\item {\tt Replacement}: This works like the {\tt Replacement} option for {\tt RandomPoints}.
\item {\tt MaxCoordinatesToReplace}: The user can specify that only a specified number of coordinates should be non-monomial (assuming {\tt Homogeneous => true}).  
\item {\tt Homogeneous}:  Setting {\tt Homogeneous => false} will cause degree zero terms to be added to modified coordinates (including monomial coordinates).		
\end{itemize} 
\begin{example}  We demonstrate some of these options.
  {{\small\color{blue}
  \begin{quote}
\begin{verbatim}
	i3 : R = ZZ/11[x, y, z];
	
	i4 : randomCoordinateChange(R)
	
           ZZ
o4 = map(R,--[x, y, z],{4x + 5y - 5z, 3x - 4y - 3z, 4x})                     
           11
           
                    ZZ
o4 : RingMap R <--- --[x, y, z]
                    11
                    
	i5 : matrix randomCoordinateChange(R, MaxCoordinatesToReplace => 1)
	
	o5 = | x -x-4y-5z y |
	
	i6 : matrix randomCoordinateChange(R, MaxCoordinatesToReplace => 1, 
	Homogeneous => false)
	
	o6 = | x-3 z-5 -x+3y-4z+2 |
\end{verbatim}
\end{quote}
    }}
  
\end{example}

\section{Other functions: {\tt genericProjection}, {\tt projectionToHypersurface}}\label{projectionfunctions}

We include two functions providing customizable projections. We describe them here. 

\subsection{genericProjection} 
This function finds a random (somewhat, depending on options) generic projection of the ring or ideal.
The typical usages are as follows: 
\vspace{0.5em}

-- {\tt genericProjection(n, I)}

-- {\tt genericProjection(n, R)}

\vspace{0.5em}
\noindent where 
$I$ is an ideal 
in a polynomial ring, 
$R$ can denote a quotient of a polynomial ring and 
$n\in \mathbb{Z}$ is
an integer specifying how many dimensions to drop.  Note that this function makes no attempt to verify that the projection is actually generic with respect to the ideal.

This gives the projection map from $\mathbb{A}^N \mapsto\mathbb{A}^{N-n}$ and the defining ideal of the projection of $V(I)$. If no integer $n$ is provided then it acts as if $n = 1$.

\begin{example}	
    We project a curve in $4$-space to one in $2$-space.
    ~~
  {{\small\color{blue}
  \begin{quote}
\begin{verbatim}
i1 : R = ZZ/5[x, y, z, w];

i2 : I = ideal(x, y^2, w^3 + x^2);

i3 : genericProjection(2, I)

            ZZ                                         2          2
o3 = (map(R,--[z, w],{- x - 2y - z, - y - 2z}), ideal(z  - z*w - w ))                   
            5
\end{verbatim}
\end{quote}
    }}
\end{example}\vspace{-1em}
\noindent
Alternatively, instead of {$I$}, we may pass it a quotient ring.  It will then return the inclusion of the generic projection ring into the given ring, followed by the source of that inclusion.  

This method works by calling {\tt randomCoordinateChange} (\Cref{helper}) before dropping variables.  It passes the options {\tt Replacement}, {\tt MaxCoordinatesToReplace}, {\tt Homogeneous} to that function.

\subsection{\tt projectionToHypersurface} This function creates a projection to a hypersurface. The typical usages are as follows: 

\vspace{0.5em}
-- {\tt projectionToHypersurface I},

-- {\tt projectionToHypersurface R} 
\vspace{0.5em}

\noindent where $I$ is an ideal in a polynomial ring, respectively, $R$ is a quotient of a polynomial ring. The output is
a list with two entries: the generic projection map and the ideal (respectively the ring).

It differs from {\tt genericProjection(codim I - 1, I)} as it only tries to find a hypersurface equation that vanishes along the projection, instead of finding one that vanishes exactly at the projection.  This can be faster and can be useful for finding points.  The same approach was used in the {\tt point} command in the {\tt Cremona} package \cite{CremonaSource}.  If we already know the codimension is {\tt c}, we can set {\tt Codimension=>c} so the function does not compute it.

\section{An application:  {\tt findANonZeroMinor}, {\tt extendIdealByNonZeroMinor}}

As mentioned in the introduction, the two functions in this section will provide further tools for computing singular locus, in addition to those available in the package {\tt FastLinAlg}. 

\subsection{\tt findANonZeroMinor:}\label{findANonZeroMinor} The typical usage of this function is as follows: 

\vspace{0.5em}
-- {\tt findANonZeroMinor(n, M, I)} 
\vspace{0.5em}

\noindent where $I$ is an ideal
in a polynomial ring over {\tt QQ} or {\tt ZZ/p} for $p$ prime, $M$ is a matrix
over the polynomial ring and $n\in \mathbb{Z}$ denotes the size of the minors of interest.

The function outputs the following:

-- randomly chosen point $P$ in $V(I)$ which it finds using {\tt randomPoints}.

-- the indexes of the columns of $M$ that stay linearly independent upon plugging $P$ into $M$, 

-- the indices of the linearly independent rows of the matrix extracted from $M$ in the above step, 

-- a random $n\times n$ sub-matrix of $M$ that has full rank at $P$.

Besides the options from {\tt randomPoints} which are automatically passed to that function, the user may also provide the following additional information: 

\begin{description}
    \item[\tt MinorPointAttempts => ZZ] 
  This controls how many points at which to check the rank.
\end{description}
\begin{example}
  We demonstrate this function.
  {{\small\color{blue}
  \begin{quote}
\begin{verbatim}
		i3 : R = ZZ/5[x, y, z];
		
		i4 : I = ideal(random(3, R) - 2, random(2, R))
		
		              3     2        2    3    2              2        2       2
		o4 = ideal (2x  - 2x y + 2x*y  + y  + x z - 2x*y*z + y z - 2x*z  + 2y*z 
		-------------------------------------------------------------------
		   3                      2
		- z  - 2, - 2x*y - x*z - z )
		
		o4 : Ideal of R
		
		i5 :  M = jacobian(I)
		
		o5 = {1} | x2+xy+2y2+2xz-2yz-2z2   -2y-z |
		     {1} | -2x2-xy-2y2-2xz+2yz+2z2 -2x   |
		     {1} | x2-2xy+y2+xz-yz+2z2     -x-2z |
		
		             3       2
		o5 : Matrix R  <--- R
		
		i6 : findANonZeroMinor(2, M, I, Strategy => GenericProjection)
		
		o6 = ({-2, 1, 1}, {0, 1}, {0, 1}, {1} | x2+xy+2y2+2xz-2yz-2z2   -2y-z |)
		{1} | -2x2-xy-2y2-2xz+2yz+2z2 -2x   |
\end{verbatim}
\end{quote}  
    }}
\end{example}

\subsection{\tt extendIdealByNonZeroMinor:}\label{extendIdealByNonZeroMinor} The typical usage is 

\vspace{0.5em}
-- {\tt extendIdealByNonZeroMinor(n, M, I)} 

\vspace{0.5em}
\noindent where $n,M,I$ are same as before. This function finds a submatrix of size $n\times n$ using \\{\tt findANonZeroMinor};  
it extracts the last entry of the output, finds its determinant and
adds it to the ideal $I$, thus extending $I$.   It has the same options as {\tt findANonZeroMinor}.

One can use this function to show that rings are $(R_1)$, that is, regular in codimension $1$.

\begin{example}
    Consider the following 3-dimensional example (which is $(R_1)$, regular in codimension 1) where computing the dimension of the singular locus takes around $30$ seconds as there are $15500$ minors of size $4 \times 4$ in the associated $7 \times 12$ Jacobian matrix.  However, we can use this function to quickly find interesting minors.
  {{\small\color{blue}
  \begin{quote}
\begin{verbatim}
		i2 : T = ZZ/101[x1, x2, x3, x4, x5, x6, x7];
		
		i3 : I =  ideal(x5*x6-x4*x7,x1*x6-x2*x7,x5^2-x1*x7,x4*x5-x2*x7,x4^2-x2*x6,x1*x4-x2*x5,
		x2*x3^3*x5+3*x2*x3^2*x7+8*x2^2*x5+3*x3*x4*x7-8*x4*x7+x6*x7,
		x1*x3^3*x5+3*x1*x3^2*x7+8*x1*x2*x5+3*x3*x5*x7-8*x5*x7+x7^2,
		x2*x3^3*x4+3*x2*x3^2*x6+8*x2^2*x4+3*x3*x4*x6-8*x4*x6+x6^2,
		x2^2*x3^3+3*x2*x3^2*x4+8*x2^3+3*x2*x3*x6-8*x2*x6+x4*x6,
		x1*x2*x3^3+3*x2*x3^2*x5+8*x1*x2^2+3*x2*x3*x7-8*x2*x7+x4*x7,
		x1^2*x3^3+3*x1*x3^2*x5+8*x1^2*x2+3*x1*x3*x7-8*x1*x7+x5*x7);
		
		o3 : Ideal of T
		
		i4 : M = jacobian I;

             7       12
o4 : Matrix T  <--- T          

i5 : i = 0; J = I;

o6 : Ideal of T

i7 : elapsedTime(while (i < 10) and dim J > 1 do (
                                i = i + 1;                       
                                J = extendIdealByNonZeroMinor(4, M, J)));
           -- 0.640164 seconds elapsed        

i8 : dim J

o8 = 1

i9 : i

o9 = 5
\end{verbatim}
\end{quote}
    }}    
    \noindent
    In this particular example, there tend to be about $5$ associated primes when adding the first minor to $J$, and so one expects about $5$ steps as each minor will typically eliminate one of those primes.
\end{example} 

There is some similar functionality computing partial Jacobian ideals obtained via heuristics (as opposed to actually finding rational points) in the package {\tt FastLinAlg}, see \cite{FastLinAlgSource}.  That package now uses the functionality contained here in {\tt RandomPoints} in some of its functions.


\section{The latest version}

The latest version of the package is available in the {\tt FastLinAlg} branch of 

\begin{center}
    {\scriptsize\tt\url{https://github.com/Macaulay2/Workshop-2020-Cleveland/blob/FastLinAlg/FastLinAlg/M2/RandomPoints.m2}}
\end{center}

\bibliographystyle{alpha}
\bibliography{MainBib}

\end{document}